\documentclass[12pt,showkeys,superscriptaddress,nofootinbib]{revtex4}
\usepackage{amsmath,amssymb,euscript,yfonts,psfrag,latexsym,dsfont,graphicx,bbm,color,amstext,wasysym,wrapfig,soul,framed,float}
\usepackage[utf8]{inputenc}
\usepackage[export]{adjustbox}

\graphicspath{{/},{figures/}}

\newcommand{\cA}{{\mathcal A}}
\newcommand{\cB}{{\mathcal B}}
\newcommand{\cE}{{\mathcal E}}

\newcommand{\cD}{{\mathcal D}}

\newcommand{\cK}{{\mathcal K}}

\newcommand{\cC}{{\mathcal C}}
\newcommand{\bbP}{\mathbb{P}}
\newcommand{\bbR}{\mathbb{R}}

\newcommand{\mD}{{\mathbb D}}
\newcommand{\bbD}{{\mD}}

\newcommand{\support}{{\rm Supp}}

\definecolor{llgrey}{rgb}{0.9,0.9,0.9}
\definecolor{lgrey}{rgb}{0.6,0.6,0.6}
\definecolor{lred}{rgb}{0.9,0.7,0.7}

\newtheorem{thm}{Theorem}

\newtheorem{lemma}{Lemma}

\newtheorem{problem}{Problem}

\newtheorem{defn}{Definition}

\begin{document}
\title{Entropy-regularized optimal transport over networks with incomplete marginals information\footnote{Supported by NYUAD grant 76/71260/ADHPG/AD356/51119}}
\author{Aayan Masood Pathan}\affiliation{Division of Science, NYUAD, U.A.E.}
\email{amp10039@nyu.edu,mp5158@nyu.edu}

\author{Michele Pavon}\affiliation{Division of Science, NYUAD, U.A.E.}

\begin{abstract} {\bf Abstract.} We study in this paper optimal mass transport over a strongly connected, directed graph on a given discrete time interval. Differently from previous literature, we do not assume full knowledge of the initial and final goods distribution over the network nodes. In spite of the meager information, we show that it is possible to characterize the most likely flow in two important cases: The first one is when the initial and/or final distribution is only known on  proper subsets of the nodes. The second case is when only some moments of the marginal distributions are known. As an important by-product, we determine the most likely initial and final marginals on the whole state space. 
\end{abstract}
\maketitle

\section{Introduction}
Problems of optimal transport over networks have been studied in the past fifteen years from several different angles, see e.g \cite{pavon2010discrete,BCCNP,chen2016robust,leoD,chen2017efficient,Pey19,Z,HCK,CLZ,DB,HRCK,KDO,LCFS,WFBR,CGKS}. 
Consider a road or computer network modeled as a directed graph where the nodes represent cities or junctions. Suppose the total number of traveling vehicles/packets is known during a discrete time interval $[0,N]$. Suppose also that at $t=0$ the vehicle/packets distribution can only be observed  at certain nodes and similarly at the final time $N$. We are interested in determining the flow on $[0,N]$.  Given that the problem is clearly {\em ill-posed}, we ask: Can we  at least determine what flow is {\em most likely}? In spite of the apparent importance, this problem has so far  received limited attention in the literature.  In the
recent paper \cite{Mascherpa2023}, authors develop
an application to estimating the pollution spreading in a
network as a discrete multi-marginal Schr\"odinger bridge
where the information on the intermediate marginals is
only on a {\em fixed subset} of the state space. Full observation of the
initial marginal is however assumed. The
dual problem is solved in \cite{Mascherpa2023} through
a block coordinate ascent method. While we study the
simpler problem with initial-final marginals, we do not
assume full knowledge of the initial marginal.  Thus, our
solution entails also the most likely completion of the initial
and final marginals. Moreover, the two marginals may be known on different subsets of the state space. The role of the initial marginal of the prior distribution on paths, which is missing in \cite{Mascherpa2023}, is clarified in this paper.  Moreover, the normalization condition of the marginal
has to be explicitly imposed in our formulation of the problem. Without that, the solution does not admit an immediate large deviation ({\em most likely}) interpretation (see, however, e.g. \cite{chen2021most}). Indeed, in the simulations of Section \ref{Example}, without the normalization condition the mass of traveling on the network does not sum to one. 
Notice that, on the subsets where the marginal distribution is known,  it is natural to impose a {\em hard constraint}. This is different from adding {\em divergences} to the cost function as done in Problem 3.1 in \cite{liero2018optimal}, and in \cite{chizat2018scaling}. 

Developing a natural modification of regularized Optimal Mass Transport (OMT) (Schr\"odinger Bridge Problem (SBP)), we show in this paper that a satisfactory answer can indeed be provided. As a by-product, the optimal flow on paths provides the most likely  {\em completion} of the initial and final marginals, a result which appears of considerable potential. The paper is outlined as follows. In Section \ref{maxent}, we collects some basic results on entropy-regularized optimal transport over networks in the case when full information of the initial or/and final marginal is available. Section \ref{partial} contains our main results. In the case where information on one or both marginals is only available on proper subsets of the state space, we characterize the optimal flow. Existence and uniqueness  of the optimal solution is established proving in Subsection \ref{ITAL} convergence of a suitable iterative scheme. To illustrate our main result, we solve in Section \ref{Example} a numerical example both when the prior distribution on paths is of the Boltzmann type and when it is the Ruelle-Bowen measure. In the following Section \ref{sec:prior}, we briefly comment on how prior transition probabilities are modified when the prior random vector is a ``time-window" of a Markov chain. In Section \ref{moments}, we outline how to solve the maximum entropy problem when only some moments of the marginals are known. Finally, some complementary, rather straightforward results/extensions are presented in Appendices \ref{mean}, \ref{variance}, \ref{continuous}.

\section{Background on discrete Schr\"odinger Bridges}\label{maxent}

Consider a  finite state space, for notational simplicity $\mathcal X=\{1,2,\ldots,n\}$. We are concerned with probability distributions $P$ on ``paths" $x=(x_0,x_1,\ldots,x_N)$ in $\mathcal X^{N+1}$.
We write  $\mathcal P$ for the simplex of all such distributions. 
For $P\in \mathcal P$, we denote by
\[
p(s,x_s;t,x_t):=P(X(t)=x_t \mid X(s)=x_s), \quad 0\le s < t\le N, \quad x_s, x_t\in{\cal X}.
\]
its transition probabilities. We also use $p$, but indexed, to denote marginals. Thus,
\[
p_{t}(x_t):=P(X(t)=x_t),
\]
and similarly, for two-time marginals,
\[
p_{st}(x_s,x_t):=P(X(s)=x_s, X(t)=x_t).
\]

This formulation  represents the case where the original phase space in  Schr\"{o}dinger's problem for Brownian particles \cite{S1,S2} has undergone a  ``coarse graining" in Boltzmann's style and where time has also been discretized. Thus, the {\em a priori} model is now given by a  distribution $P\in\mathcal P$ and suppose that in experiments an initial  marginal $\nu_0$ or final marginal $\nu_T$ or both, respectively, have been observed that differ from the marginals $ p_0$ and $p_T$ of the prior distribution $P$. We denote by ${\mathcal P}(\nu_0)$, ${\mathcal P}(\nu_T)$, ${\mathcal P}(\nu_0,\nu_T)$, respectively, the subsets of $\mathcal P$ having the prescribed marginals. We seek a distribution in one of these subsets which is close to the given prior $P$. Large deviation reasoning \cite{sanov1957largedeviations} requires that we employ as ``distance'' the {\em relative entropy}:

\begin{defn}{\em Let $P,Q\in \mathcal P$ the simplex  of probability distributions on ${\cal X}^{N+1}$, and let $x=(x_0,x_1,\ldots,x_T)$. If $P(x)=0\Rightarrow Q(x)=0$, we say that the {\em support} of $Q$ is contained in the support of $P$ and write \[\support (Q)\subseteq \support (P).\]
The {\em Relative Entropy} of $Q$ from $P$ is defined to be
\begin{equation}\label{KLdist}\bbD(Q\|P)=\left\{\begin{array}{ll} \sum_{x\in{\cal X}^{T+1}}Q(x)\log\frac{Q(x)}{P(x)}, & \support (Q)\subseteq \support (P),\\
+\infty , & \support (Q)\not\subseteq \support (P).\end{array}\right.,
\end{equation} 
where, by definition,  $0\cdot\log 0=0$.}
\end{defn}

The relative entropy is also known as the {\em divergence} or {\em Kullback-Leibler index}. {As is well known, $\bbD(Q\|P)\ge 0$ and $\bbD(Q\|P)=0$ if and only if $Q=P$. }Given this notion of distance,
we seek
a probability law $Q^*$ in a suitable family 
which is closest to the prior distribution $P$ in this sense. 

\subsection{A decomposition of relative entropy}
For $P,Q\in \mathcal P$, let 
$Q_{x_0,x_N}(\cdot)=Q(\cdot|X_0=x_0, X_N=x_N)$ and similarly for $P$. Then, for both distributions, we have
$$Q(x_0,x_1,\ldots,x_N)=Q_{x_0,x_N}(x_0,x_1,\ldots,x_N)q_{0N}(x_0,x_N),
$$
where we have assumed that $q_{0N}$ (resp. $p_{0N}$)  is everywhere positive on $\mathcal X\times \mathcal X$.
We get
\begin{equation}\label{DECO}
\bbD(Q\|P)=\sum_{x_0x_N}q_{0N}(x_0,x_N)\log \frac{q_{0N}(x_0,x_N)}{p_{0N}(x_0,x_N)}+\sum_{x\in{\cal X}^{N+1}}Q_{x_0,x_N}(x)\log \frac{Q_{x_0,x_N}(x)}{P_{x_0,x_N}(x)} q_{0N}(x_0,x_N).
\end{equation}
This is the sum of two nonnegative quantities. The second becomes zero if and only if $Q_{x_0,x_N}(x)=P_{x_0,x_N}(x)$ for all $x\in{\cal X}^{N+1}$.  Thus, in any minimization of $\bbD(Q\|P)$ with constrained initial or final or both marginals, $Q^*_{x_0,x_N}(x)=P_{x_0,x_N}(x)$. In particular, this shows that when both marginals are imposed and they are   delta distributions, the solution is simply obtained from the prior through {\em conditioning}. From this point of view, Schr\"{o}dinger bridges with general (initial, final or both) marginals appear as a sort of  ``soft conditioning" of the prior thereby generalizing the Bayesian approach.

\subsection{Half-bridges}
When only one marginal is prescribed, we use the terminology {\em half-bridge} problem. Their solution is straightforward \cite[Subsection IIIB]{pavon2010discrete}. Nevertheless, we outline the argument for the purpose of comparison with the incomplete information case treated in Subsection \ref{hbpi}. 

Using $p_{0N}(x_0,x_T)=p_0(x_0)p(0,x_0;N,x_N)$ for both $P$ and $Q$, we immediately get the decomposition
\[\bbD(q_{0N}\|p_{0N})=\bbD(q_0\|p_0)+\sum_{x_0}\bbD(q(0,x_0;N,x_N)\|p(0,x_0;N,X_N))q_0(x_0).
\]
Let $\mathcal P(\nu_0)$ be the family of distributions in $\mathcal P$ having initial marginal $\nu_0$. Then, considering that both terms in the above decomposition are non negative, the minimizer of $\bbD(Q\|P)$ over $\mathcal P(\nu_0)$ is the distribution with initial marginal $\nu_0$ and $q^*(0,x_0;N,x_N)=p(0,x_0;N,x_N)$. Similarly, let $\bar{p}(N,x_N;0,x_0)=P(X(0)=x_0 \mid X(N)=x_N),$ be the reverse-time transition probability. Then, using,  $p_{0N}(x_0,x_N)=p_N(x_N)\bar{p}(N,x_N;0,x_0)$ for both $P$ and $Q$, we get 
\[\bbD(q_{0N}\|p_{0N})=\bbD(q_N\|p_N)+\sum_{x_0}\bbD(\bar{q}(N,x_N;0,x_0)\|\bar{p}(N,x_N;0,x_0))q_N(x_N).
\]
Let $\mathcal P(\nu_N)$ be the family of distributions in $\mathcal P$ having final marginal $\nu_N$. Then, considering that both terms in the above decomposition are non negative, the minimizer of $\bbD(q_{0N}\|p_{0N})$ over $\mathcal P(\nu_N)$ is the distribution with final marginal $\nu_N$ and $\bar{q}^*(N,x_N;0,x_0)=\bar{p}(N,x_N;0,x_0)$. It is interesting to compute  $q^*(0,x_0;N,x_N)$. We have
\[q^*_0(x_0)q^*(0,x_0;N,x_N)=q^*_{0N}(x_0,x_N)=\nu_N(x_N)\bar{q}^*(N,x_N;0,x_0)=\nu_N(x_N)\bar{p}(N,x_N;0,x_0),
\]
which gives (assuming that all the one-time marginals are everywhere positive)
\begin{equation}\label{FBR1}
q^*(0,x_0;N,x_N)=\frac{\nu_N(x_N)}{q^*_0(x_0)}\bar{p}(N,x_N;0,x_0)=\frac{p_0(x_0)}{q^*_0(x_0)}\frac{\nu_N(x_N)}{p_N(x_N)}p(0,x_0;N,x_N).
\end{equation}
Define
\[\varphi(t,x_t)=\frac{q^*_t(x_t)}{p_t(x_t)}, \quad t=0,1,\ldots,N.
\]
Then, (\ref{FBR1}) can be rewritten as
\begin{equation}\label{FBR2}
q^*(0,x_0;N,x_N)=\frac{\varphi(N,x_N)}{\varphi(0,x_0)}p(0,x_0;N,x_N).
\end{equation}
Finally, observe that the above gives
\begin{equation}\label{phi}\sum_{x_N}p(0,x_0;N,x_N)\varphi(N,x_N) = \varphi(0,x_0)\sum_{x_N}q^*(0,x_0;N,x_N)=\varphi(0,x_0).
\end{equation}

\subsection{General Schr\"odinger Bridge}
\begin{problem}\label{prob:optimization}{\em Assume that $p(0,\cdot\,; T,\cdot)$ is everywhere positive on $\cal X\times \cal X$. Determine
 \begin{eqnarray}\label{eq:optimization}
Q^*={\rm argmin}\{ \bbD(Q\|P)&\mid& Q\in {\mathcal P}(\nu_0,\nu_N)
\}.\nonumber
\end{eqnarray}}
\end{problem}
It turns out that if there is at least one $Q$ in
$\mathcal P(\nu_0,\nu_T)$ such that
$\bbD(Q\|P)<\infty$,
there exists a unique minimizer $Q^*$  called
{\em the Schr\"{o}dinger bridge} from $\nu_0$ to $\nu_N$ over $P$. The latter can be characterized as follows (see e.g. \cite[Subsection IIIA]{GP}):
\begin{equation}\label{NewTransition}q^*(0,x_0;N,x_N)=\frac{\varphi(N,x_N)}{\varphi(0,x_0)}p(0,x_0;N,x_N)
\end{equation}
where $\varphi$ and $\hat{\varphi}$ {must satisfy}
\begin{eqnarray}\label{OSchonestep1}
\hat{\varphi}(N,x_N)=\sum_{x_0}p(0,x_0;N,x_N)\hat{\varphi}(0,x_0),\\\label{OSchonestep2}\quad \varphi (0,x_0):=\sum_{x_N}p(0,x_0;N,x_N)\varphi(N,x_N)
\end{eqnarray}
with the boundary conditions
\begin{equation}\label{OBConestep}
\varphi(0,x_0)\cdot\hat{\varphi}(0,x_0)=\nu_0(x_0),\quad \varphi(N,x_N)\cdot\hat{\varphi}(N,x_N)=\nu_N(x_N),\quad \forall x_0, x_N\in\mathcal X.
\end{equation}
Notice that the normalization of $q^*_{0T}(x_0,x_N)=\nu_0(x_0)q^*(0,x_0;N,x_N)$ follows from the fact that the marginals are probability distributions. Indeed,
\[\sum_{x_0,x_N}q^*_{0T}(x_0,x_N)=\sum_{x_0}\nu_0(x_0)=1.
\]
The question of existence and uniqueness\footnote{Uniqueness is always up to multiplying $\varphi$ by a positive constant $c$ and dividing $\hat{\varphi}$ by the same constant.}  of  functions $\hat{\varphi}$, $\varphi$ satisfying (\ref{OSchonestep1})-(\ref{OSchonestep2})-(\ref{OBConestep}) can be established showing contractivity in the Hilbert metric of an iterative scheme, see e.g \cite[Section IIID]{GP}.

We can now consider a simple example with a Boltzmann prior:
 \begin{equation}\label{Bolt}P_B(x_0,x_1,\ldots,x_N)=Z(T)^{-1}\prod_{t=0}^{N-1}\exp\left[-\frac{l_{x_tx_{t+1}}}{T}\right].
\end{equation}
Here, for simplicity, the length $l_{x_tx_{t+1}}$ of the edge is always $1$ if the edge exists and $+\infty$ otherwise. As is well known, as $T\searrow 0$, the Boltzmann distribution tends to concentrate on the absolute minima of the Hamiltonian, here the shortest paths. We choose here $T=1$.
The topology of the directed graph  is as follows:

\begin{figure}[H]
    \begin{center}
    \includegraphics[width=10cm]{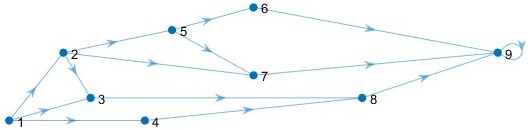}
    \end{center}
    \caption{Graph topology}
    \label{fig:topology_complete}
\end{figure}

We consider a source-sink problem with  final time $N=4$. Here, node $1$ is the source ($\nu_0(x_0)=\delta_1(x_0)$) and node $9$ is the sink ($\nu_4(x_4)=\delta_9(x_4)$).  Because of the loop on edge $9$, goods may reach the sink in less than three time units and then remain there. The below simulation shows the optimal mass evolution over time (from top to bottom).

\begin{figure}[H]
    \begin{center}
    \includegraphics[width=15cm]{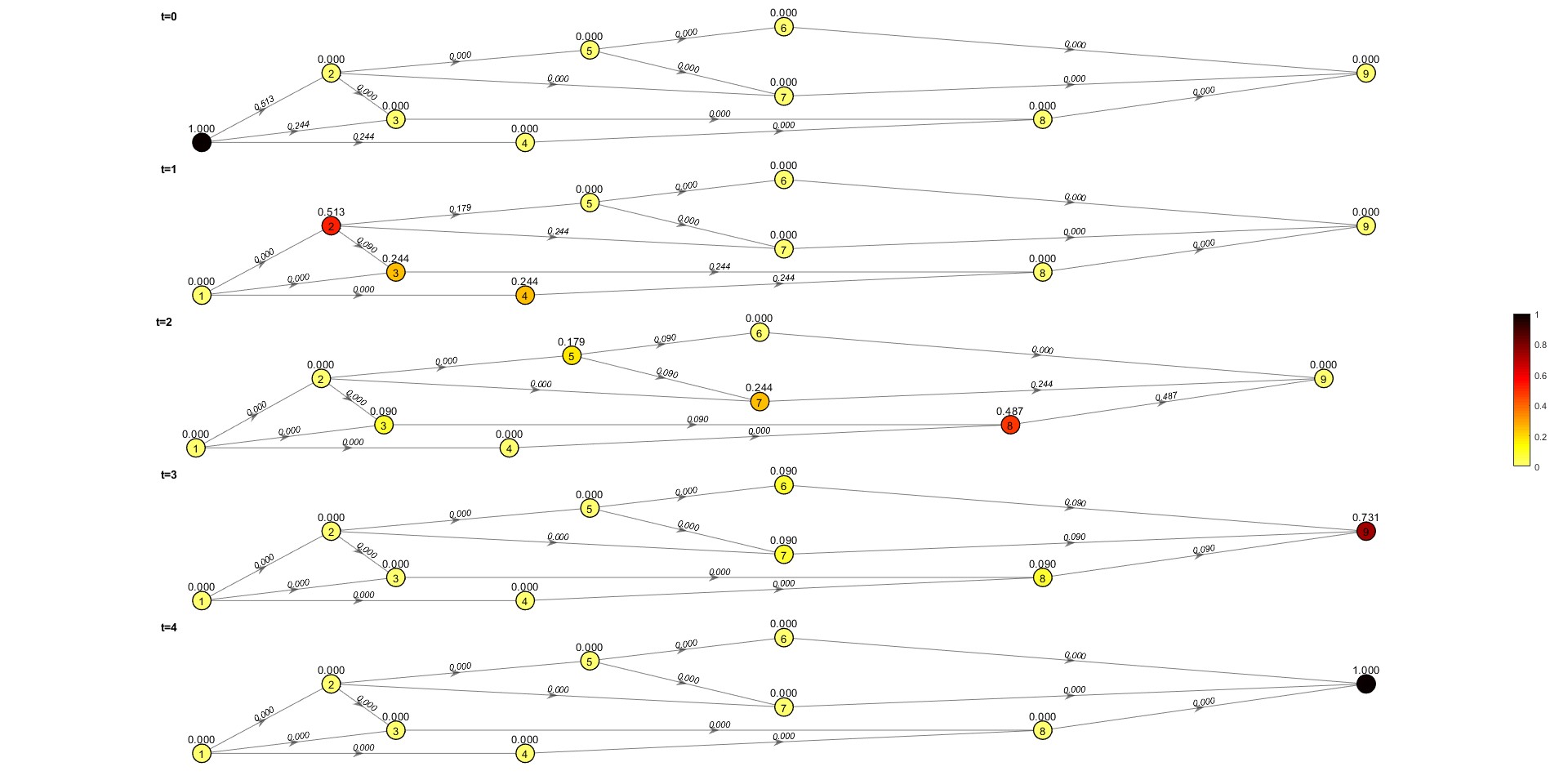}
    \end{center}
    \caption{Optimal network flow in the source-sink problem with $t=0,1,2,3,4$ and Boltzmann prior.}
    \label{fig:complete}
\end{figure}

\noindent
We notice that the paths from source  to  sink take 3 or 4 steps, and that most of the mass, due to the low value of $T$,  travels on the shortest paths.

\section{Partial knowledge of the marginals}\label{partial}

We consider now the rather common and important situation where the initial and/or final distributions  $\nu_0$ and $\nu_N$, are only known on $\mathcal X_0\subseteq \mathcal X$ and $\mathcal X_N\subseteq \mathcal X$, respectively. Let $\mathcal P_{0N}$ be the family of probability distributions over $\mathcal X\times\mathcal X$. As for the case of full marginal knowledge, we consider first the half-bridge problem.

\subsection{Half-bridges with partial information}\label{hbpi}
We consider first the situation where the initial marginal is only known on  $\mathcal X_0\subset \mathcal X$. More specifically, let $\rho_0$ be a positive measure on  $\mathcal X_0$ satisfying 
\begin{equation}\label{0cond}0<\alpha=\sum_{x\in \mathcal X_0}\rho_0(x)<1.
\end{equation}
Clearly, both situations $\alpha =0$ and $\alpha=1$ lead to a trivial solution of the maximum entropy problem.  For any set $A$, let $\mathds{1}_A$ denote the indicator function defined
\[\mathds{1}_A(x)=\left\{\begin{array}{ll} 1,& {\rm if}\;x\in A\\ 0, & {\rm if}\;x\not\in A.\end{array}\right.
\]Let $\mathcal P(\rho_0)$ be the family of distributions in $\mathcal P$ whose initial marginal $q_0$ coincides with $\rho_0$ on $\mathcal X_0$, i.e. $\left[q_0(x)-\rho_0(x)\right]\mathds{1}_
{\mathcal X_0}\equiv 0$. By the same argument as with full initial marginal knowledge, the minimizer of $\bbD(Q\|P)$ over $\mathcal P(\rho_0)$ is the distribution with $q^*(0,x_0;T,x_T)=p(0,x_0;T,x_T)$ and minimizing $\bbD(q_0\|p_0)$. Since
\[\bbD(q_0\|p_0)=\sum_{x_0\in\mathcal X_0}\log\frac{q_0(x_0)}{p_0(x_0)}q_0(x_0) + \sum_{x_0\in\mathcal (X_0)^c}\log\frac{q_0(x_0)}{p_0(x_0)}q_0(x_0)
\]
and the first term is constant over $\mathcal P(\rho_0)$, the problem reduces to minimize $\sum_{x_0\in\mathcal (X_0)^c}\log\frac{q_0(x_0)}{p_0(x_0)}q_0(x_0)$ over $q_0$ coinciding with $\rho_0$ on $\mathcal X_0$ subject to the constraint 
\[\sum_{x_0\in(\mathcal X_0)^c} q_0(x_0)=1-\alpha.
\]
The Lagrangian is  
\[\mathcal L(q_0;\lambda) = \sum_{x_0\in\mathcal (X_0)^c}\log\frac{q_0(x_0)}{p_0(x_0)}q_0(x_0)+\lambda(\sum_{x_0\in(\mathcal X_0)^c} q_0(x_0)-1+\alpha).
\]Setting the first variation with respect to $q_0$ equal to zero, we get for $x_0\in (\mathcal X_0)^c$, $q_0^*(x_0)=p_0(x_0)\exp[-\lambda-1]$.Hence, on $(\mathcal X_0)^c$, $q_0^*(x_0)=p_0(x_0)c_0$, where the constant $c_0$ is given by 
\[c_0=\frac{1-\alpha}{\sum_{x_0\in(\mathcal X_0)^c} p_0(x_0)},
\]
since both $p_0$ and $q_0^*$ sum to one over $\mathcal X$. In conclusion, $q^*(0,x_0;T,x_T)=p(0,x_0;T,x_T)$  and
\[q_0^*(x_0) = \left\{\begin{array}{ll} \rho_0(x_0),& {\rm if}\;x_0\in \mathcal X_0\\ p_0(x_0)c_0, c_0=\frac{1-\alpha}{\sum_{x_0\in(\mathcal X_0)^c} p_0(x_0)}& {\rm if}\;x_0\in (\mathcal X_0)^c.\end{array}\right.
\]
We consider now the situation where the final marginal is only known on  $\mathcal X_N\subset \mathcal X$. More specifically, let $\rho_N$ be a positive measure on  $\mathcal X_N$ satisfying 
\begin{equation}\label{Tcond}0<\beta=\sum_{x\in \mathcal X_N}\rho_N(x)<1.
\end{equation}
Let $\mathcal P(\rho_N)$ be the family of distributions in $\mathcal P$ whose final marginal $q_N$ coincides with $\rho_N$ on $\mathcal X_N$, i.e. $\left[q_N(x)-\rho_N(x)\right]\mathds{1}_
{\mathcal X_N}\equiv 0$. By the same argument as with full final marginal knowledge, the minimizer of $\bbD(Q\|P)$ over $\mathcal P(\rho_N)$ is the distribution with $\bar{q}^*(N,x_N;0,x_0)=\bar{p}(N,x_N;0,x_0)$ and minimizing $\bbD(q_N\|p_N)$. By the same argument as before, the problem reduces to minimize $\sum_{x_N\in\mathcal (X_N)^c}\log\frac{q_N(x_N)}{p_N(x_N)}q_N(x_N)$ over $q_N$ coinciding with $\rho_N$ on $\mathcal X_N$ subject to the constraint 
\[\sum_{x_N\in(\mathcal X_N)^c} q_N(x_N)=1-\beta.
\]
As before, we get $\bar{q}^*(N,x_N;0,x_0)=\bar{p}(N,x_N;0,x_0)$   and
\[q_N^*(x_N) = \left\{\begin{array}{ll} \rho_N(x_N),& {\rm if}\;x_N\in \mathcal X_N\\ p_N(x_N)c_N, c_N=\frac{1-\beta}{\sum_{x_N\in(\mathcal X_N)^c} p_N(x_N)}& {\rm if}\;x_N\in (\mathcal X_N)^c.\end{array}\right.
\]
As in the case of full marginal information, we compute $q^*(0,x_0;T,x_T)$. We have
\[q^*_0(x_0)q^*(0,x_0;N,x_N)=q^*_{0N}(x_0,x_N)=q^*_N(x_N)\bar{q}^*(N,x_N;0,x_0)=q^*_N(x_N)\bar{p}(N,x_N;0,x_0),
\]
which gives (assuming that all the one-time marginals are everywhere positive)
\begin{equation}\label{FBR}
q^*(0,x_0;N,x_N)=\frac{q^*_N(x_N)}{q^*_0(x_0)}\bar{p}(N,x_N;0,x_0)=\frac{p_0(x_0)}{q^*_0(x_0)}\frac{q^*_N(x_N)}{p_N(x_N)}p(0,x_0;N,x_N).
\end{equation}
Define
\[\varphi(t,x_t)=\frac{q^*(t,x_t)}{p(t,x_t)}, \quad t=0,1,\ldots,N.
\]
Observe that
\begin{equation}\label{endvarphi}\varphi(N,x_N)= \left\{\begin{array}{ll} \frac{\rho_N(x_N)}{p_N(x_N)},& {\rm if}\;x_N\in \mathcal X_N\\ c_N, c_N=\frac{1-\beta}{\sum_{x_N\in(\mathcal X_N)^c} p_N(x_N)}& {\rm if}\;x_N\in (\mathcal X_N)^c.\end{array}\right.
\end{equation}
Then, (\ref{FBR}) can be rewritten as
\begin{equation}\label{FBR2}
q^*(0,x_0;N,x_N)=\frac{\varphi(N,x_N)}{\varphi(0,x_0)}p(0,x_0;N,x_N).
\end{equation}
Finally, observe that the above gives
\[\sum_{x_N}p(0,x_0;N,x_N)\varphi(N,x_N) = \varphi(0,x_0)\sum_{x_N}q^*(0,x_0;N,x_N)=\varphi(0,x_0).
\]

\subsection{The full bridge with partial information}\label{FBPI}
Consider now the situation where the initial marginal must coincide with $\rho_0$ on $\mathcal X_0$ and the final marginal must coincide with $\rho_N$ on $\mathcal X_N$. We call this problem the {\em Incomplete Marginal Schr\"odinger Bridge Problem} (IMSBP). As already observed, because of the decomposition (\ref{DECO}), the maximum entropy problem  becomes minimizing
\begin{equation}\label{IMSBP1}\bbD(q_{0N}\|p_{0N})=\sum_{x_0x_N}q_{0N}(x_0,x_N)\log \frac{q_{0N}(x_0,x_N)}{p_{0N}(x_0,x_N)}
\end{equation}
with respect to $q_{0N}\in\mathcal P_{0N}$ subject to the (linear) constraints
\begin{eqnarray}\label{IMSBP2}\sum_{x_N}q_{0N}(x_0,x_N)&=&\rho_0(x_0),\quad x_0\in{\mathcal X_0},
\\
\label{IMSBP3}\sum_{x_0}q_{0N}(x_0,x_N)&=&\rho_N(x_N),\quad x_N\in{\mathcal X_N},\\\label{IMSBP4}
\sum_{x_0x_N}q_{0N}(x_0,x_N)&=&1.
\end{eqnarray}

The Lagrangian function has the form
\begin{eqnarray}\nonumber&&{\cal L}(q_{0N})=\sum_{x_0x_N}q_{0N}(x_0,x_N)\log \frac{q_{0N}(x_0,x_N)}{p_{0N}(x_0,x_N)}+\sum_{x_0\in\mathcal X}\mathds{1}_{\mathcal X_0}(x_0)\lambda(x_0)\left[\sum_{x_N}q_{0N}(x_0,x_N)-\rho_0(x_0)\right]\\&&+\sum_{x_N\in\mathcal X}\mathds{1}_{\mathcal X_N}(x_N)\mu(x_N)\left[\sum_{x_0}q_{0N}(x_0,x_N)-\mathbf \rho_N(x_N)\right]+\theta\left[\sum_{x_0x_N}q_{0N}(x_0,x_N)-1\right]\label{lagrangian}
\end{eqnarray}
Setting the first variation equal to zero, we get the (sufficient) optimality conditions
\[1+\log q_{0N}^*(x_0,x_N)-\log p_{0N}(x_0,x_N)+\mathds{1}_{\mathcal X_0}(x_0)\lambda(x_0)+\mathds{1}_{\mathcal X_N}(x_N)\mu(x_N)+\theta=0, (x_0,x_N)\in (\mathcal X\times\mathcal X).
\]
Using $p_{0N}(x_0,x_N)=p_0(x_0)p(0,x_0;N,x_N)$, we get
\[\frac{q_{0N}^*(x_0,x_N)}{p(0,x_0; N,x_N)}=\left\{p_0(x_0)\exp[-1-\theta-\mathds{1}_{\mathcal X_0}(x_0)\lambda(x_0)]\right\}\left\{\exp\left[-\mathds{1}_{\mathcal X_N}(x_N)\mu(x_N)\right]\right\}=\hat{\varphi}(x_0)\cdot\varphi(x_N).
\]
Hence,  the ratio $q_{0N}^*(x_0,x_N)/p(0,x_0; N,x_N)$ factors into a function of $x_0$ times a function of $x_N$,  denoted $\hat{\varphi}(x_0)$ and $\varphi(x_N)$, respectively. More explicitly,
\begin{equation}\nonumber\hat{\varphi}(x_0)= \left\{\begin{array}{ll} p_0(x_0)\exp [-1-\theta-\lambda(x_0)],& {\rm if}\;x_0\in \mathcal X_0\\ p_0(x_0)\exp [-1-\theta],& {\rm if}\;x_0\in (\mathcal X_0)^c\end{array}\right..
\end{equation}
Moreover,
\begin{equation}\nonumber\varphi(x_N)= \left\{\begin{array}{ll} \exp \left[-\mu(x_N)\right],& {\rm if}\;x_N\in \mathcal X_N\\ 1,& {\rm if}\;x_N\in (\mathcal X_N)^c\end{array}\right.
\end{equation}
We can then write the optimal $q_{0N}^*(\cdot,\cdot)$ in  the form 
\begin{equation}\label{optimaljoint} q_{0N}^*(x_0,x_N)=\hat{\varphi}(x_0) p(0,x_0;N,x_N)\varphi(x_N),
\end{equation} 
where $\varphi$ and $\hat{\varphi}$ {must satisfy}
\begin{eqnarray}\hat{\varphi}(x_0)\sum_{x_N}p(0,x_0;N,x_N)\varphi(x_N)&=&\rho_0(x_0),\quad x_0\in\mathcal X_0,\label{opt1}\\
\varphi(x_N)\sum_{x_0}p(0,x_0;N,x_N)\hat{\varphi}(x_0)&=&\rho_N(x_N), \quad x_N\in\mathcal X_N.\label{opt2}
\end{eqnarray}
Let us define $\hat{\varphi}(0,x_0)=\hat{\varphi}(x_0)$, $\quad \varphi(N,x_N)=\varphi(x_N)$ and  
$$\hat{\varphi}(N,x_N)=\sum_{x_0}p(0,x_0;N,x_N)\hat{\varphi}(0,x_0),\quad \varphi (0,x_0):=\sum_{x_N}p(0,x_0;N,x_N)\varphi(N,x_N).
$$
Then, (\ref{opt1})-(\ref{opt2}) can be replaced by the system
\begin{eqnarray}\label{Schonestep1}
\hat{\varphi}(N,x_N)=\sum_{x_0}p(0,x_0;N,x_N)\hat{\varphi}(0,x_0),\\\label{Schonestep2}\quad \varphi (0,x_0):=\sum_{x_N}p(0,x_0;N,x_N)\varphi(N,x_N)
\end{eqnarray}
with the boundary conditions
\begin{equation}\label{BConestep}
\varphi(0,x_0)\cdot\hat{\varphi}(0,x_0)=\rho_0(x_0),\quad \forall x_0\in\mathcal X_0\quad \varphi(N,x_N)\cdot\hat{\varphi}(N,x_N)=\rho_N(x_N),\quad \forall  x_N\in\mathcal X_N.
\end{equation}
Differently from the full info on marginals case, we need here to impose the normalization of $q^*_{0N}$ in one of the two equivalent forms
\begin{equation}\label{normalization}\sum_{x_0}\hat{\varphi}(0,x_0)\varphi(0,x_0)=1, \quad {\rm or}\; \sum_{x_N}\hat{\varphi}(N,x_N)\varphi(N,x_N)=1.
\end{equation}
Existence for the system (\ref{Schonestep1})-(\ref{Schonestep2})-(\ref{BConestep})-(\ref{normalization}) is established in the next subsection by proving convergence of a suitable iterative scheme.

\subsection{An iterative algorithm}\label{ITAL}
Consider the following four maps
\begin{eqnarray}\hat\varphi(0,x_0)  \overset{\cE^\dagger}{\longrightarrow} \hat\varphi(N,x_N),\quad \hat\varphi(N,x_N) &=& \sum_{x_0}p(0,x_0;N,x_N)\hat{\varphi}(0,x_0)\\
 \hat\varphi(N,x_N)\overset{\cD_{T}}{\longrightarrow}\varphi(N,x_N),\quad \varphi(N,x_N)&=& \left\{\begin{array}{ll} \frac{\rho_N(x_N)}{\hat{\varphi}(N,x_N)},& {\rm if}\;x_N\in \mathcal X_N\\ 1,& {\rm if}\;x_N\in (\mathcal X_N)^c\end{array}\right.\\
\varphi(N,x_N)
\overset{\cE}{\longrightarrow} \varphi(0,x_0),\quad \varphi (0,x_0)&=&\sum_{x_N}p(0,x_0;N,x_N)\varphi(N,x_N)\\
 \varphi(0,x_0)\overset{\cD_{0}}{\longrightarrow}\left(\hat\varphi(0,x_0)\right)_{\rm next},\quad\hat{\varphi}(0,x_0)^{{\rm next}}&=& \left\{\begin{array}{ll} \frac{\rho_0(x_0)}{\varphi(0,x_0)},& {\rm if}\;x_0\in \mathcal X_0\\ c_0p_0(x_0),& {\rm if}\;x_0\in (\mathcal X_0)^c\end{array}\right.
\end{eqnarray}
Here, the constant $c_0$ must be such that $\sum_{x_0}\hat{\varphi}(0,x_0)^{\rm next}\varphi(0,x_0)=1$ so as to satisfy (\ref{normalization}). Let $\sum_{x_0\in\mathcal X_0}\rho_0(x_0)=\alpha$. Since $\varphi(0,x_0)\cdot\hat{\varphi}(0,x_0)=\rho_0(x_0)$ for $x_0\in\mathcal X_0$, we must have
\[\sum_{x_0\in(\mathcal X_0)^c}\hat{\varphi}(0,x_0)^{\rm next}\varphi(0,x_0)=1-\alpha.
\]
We thus find the value for $c_0$, namely

\[c_0=\frac{1-\alpha}{\sum_{x_0\in(\mathcal X_0)^c}p_0(x_0)\varphi(0,x_0)}.
\]
\noindent
Let us introduce Hilbert's projective metric. Let $\mathcal S$ be a real Banach space and let $\mathcal K$ be a closed solid cone in $\mathcal S$, i.e., $\mathcal K$ is closed with nonempty interior and is such that $\mathcal K+\mathcal K\subseteq \cK$, $\mathcal K\cap -\mathcal K=\{0\}$ as well as $\lambda \mathcal K\subseteq \mathcal K$ for all $\lambda\geq 0$. Define the partial order
\[
x\preceq y \Leftrightarrow y-x\in\mathcal K,
\]
and for $x,y\in\mathcal K\backslash \{0\}$, define
\begin{eqnarray*}
M(x,y)&:=&\inf\, \{\lambda\,\mid x\preceq \lambda y\}\\
m(x,y)&:=&\sup \{\lambda \mid \lambda y\preceq x \}.
\end{eqnarray*}
Then, the Hilbert metric is defined on $\mathcal K\backslash\{0\}$ by
\[
d_H(x,y):=\log\left(\frac{M(x,y)}{m(x,y)}\right).
\]
It is a {\em projective} metric since it remains invariant under scaling by positive constants, i.e.,
$d_H(x,y)=d_H(\lambda x,y)=d_H(x,\lambda y)$ for any $\lambda>0$ and, thus, it actually measures distance between rays and not elements.

Next result extends to our more general setting \cite[Lemma 1]{GP}.
\begin{lemma}\label{contractiveness}Assume that $p(0,\cdot; N,\cdot)$ is everywhere positive on $\mathcal X\times \mathcal X$. Then the composition of the above defined four maps
\begin{equation}\label{eq:composition}
\hat\varphi(0,x_0)  \overset{\cE^\dagger}{\longrightarrow} \hat\varphi(N,x_N)\overset{\cD_{N}}{\longrightarrow}\varphi(N,x_N)
\overset{\cE}{\longrightarrow} \varphi(0,x_0)
\overset{\cD_{0}}{\longrightarrow}\left(\hat\varphi(0,x_0)\right)_{\rm next}.
\end{equation}
is contractive in the Hilbert  metric. 
\end{lemma}
{\bf Proof} Maps $\cE$ and $\cE^\dagger$ are the same as in \cite[Lemma 1]{GP}. Thus,  $\|\cE\|_H<1$ and $\|\cE^\dagger\|_H\le 1$. Consider now $\cD_N$. By the same argument as in  \cite[Lemma 1]{GP}, the map is isometric or even contractive on $\mathcal X_N$ with respect to the Hilbert metric. Since it is {\em constant} on $(\mathcal X_N)^c$, it follows that $\|\cD_N\|_H\le 1$. Similar considerations can be made for $\cD_0$ since for $x_0\in (\mathcal X_0)^c$,  functions get mapped to the same quantity $c_0p_0(x_0)$. We conclude that $\|\cD_0\|_H\le 1$. The conclusion now follows from the elementary fact
\[
\|{ \cD}_0\circ \cE\circ \cD_N\circ\cE^\dagger\|_H\leq \|{ \cD}_0\|_H\cdot\|\cE\|_H\cdot\|\cD_N\|_H\cdot\|\cE^\dagger\|_H<1,
\]
where $\circ$ denotes composition. 
{\bf Q.E.D.}

\noindent
We are now ready to establish existence and uniqueness for the Schr\"odinger system (\ref{Schonestep1})-(\ref{Schonestep2}) with boundary conditions (\ref{BConestep}) and the normalization condition (\ref{normalization}).

\begin{thm} \label{fundtheorem} Let $\rho_0$ and $\rho_N$ be positive measures on $\mathcal X_0$ and $\mathcal X_N$, respectively, satisfying the conditions of Subsection \ref{hbpi}.  Assume that $p(0,\cdot; N,\cdot)$ is everywhere positive on $\mathcal X\times \mathcal X$. Then, there exist four vectors
$
\varphi(0,x_0),\,\varphi(N,x_N),\,\hat\varphi(0,x_0),\,\hat\varphi(N,x_N)$, indexed by $x_0, x_N\in\mathcal X,
$
with nonnegative  entries satisfying (\ref{Schonestep1})-(\ref{Schonestep2})-(\ref{BConestep})-(\ref{normalization}). The four vectors are unique up to multiplication of $\varphi(0,x_0)$ and $\varphi(N,x_N)$ by the same positive constant and division of $\hat\varphi(0,x_0)$ and $\hat\varphi(N,x_N)$ by the same constant. The solution to the (IMSBP) (\ref{IMSBP1})-(\ref{IMSBP4}) is then given by the joint initial-final distribution $q_{0,N}^*(x_0,x_N)$ in (\ref{optimaljoint}).
\end{thm}
{\bf Proof} In view of Lemma \ref{contractiveness}, $\cC={ \cD}_0\circ \cE\circ \cD_N\circ\cE^\dagger$ is contractive in the Hilbert metric. Hence, there exists a unique positive $\hat\varphi(0,\cdot)=[\hat\varphi(0,x_0)]$ so that the corresponding ray is invariant under $\cC$. That is, in the notation of \eqref{eq:composition},
\begin{eqnarray*}
(\hat\varphi(0,\cdot))_{\rm next}&=&\cC(\hat\varphi(0,\cdot))\\
&=& \lambda \hat\varphi(0,\cdot).
\end{eqnarray*}
In view of the normalization condition (\ref{normalization}), we can then use the same argument as in the proof of  \cite[Theorem 3]{GP} to conclude that $\lambda=1$.
{\bf Q.E.D.}

\section{A numerical example}\label{Example}
\subsection{Boltzmann prior}
\begin{figure}[H]
    \begin{center}
    \includegraphics[width=10cm]{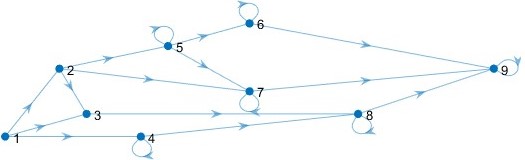}
    \end{center}
    \caption{Graph topology}
    \label{fig:topology_incomplete1}
\end{figure}

 We consider some IMSBPs when $N=4$. The graph topology is illustrated in  Figure:\ref{fig:topology_incomplete1}.

\begin{figure}
    \begin{center}
    \includegraphics[width=15cm]{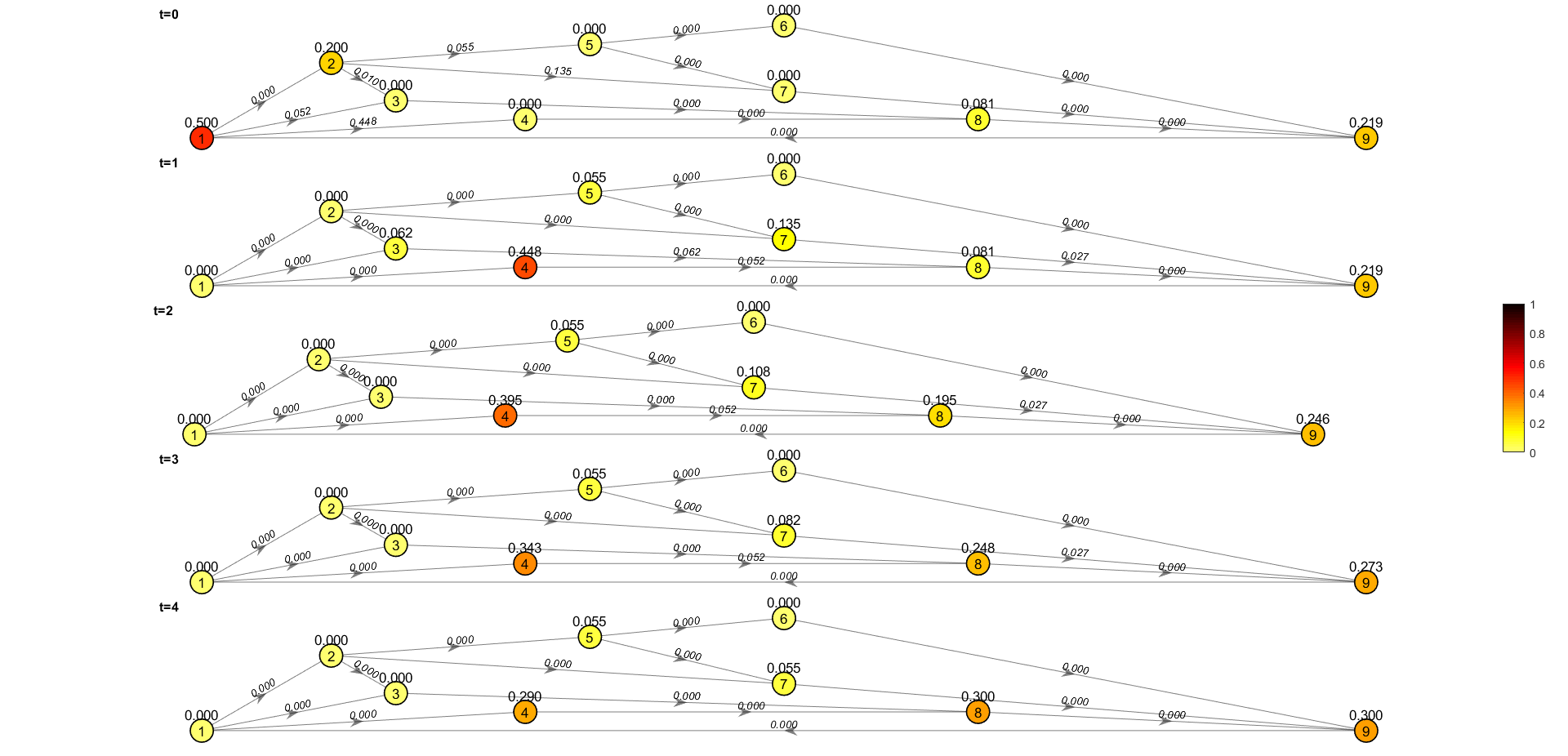}
    \end{center}
    \caption{Solution to the Schr\"odinger bridge problem with a Boltzmann prior T=0.01.}
    \label{fig:boltzmann_T_0.01}
\end{figure}

Of the marginals, we only  know $\rho_0(1)= 0.5, \rho_0(2)= 0.2$, ($\mathcal X_0=\{1,2\}$) and $\rho_4(8)= \rho_4(9)=0.3$ ($\mathcal X_4=\{8,9\}$). We consider the evolution of the distribution over the nodes when $t=0,1,2,3,4$. We first choose a Boltzmann prior (\ref{Bolt}). Choosing a low temperature, we expect mass to travel mostly on  the shortest  paths. In the following matrix, the five rows of the matrix show the mass distribution at times $t=0, 1, 2 ,3,4$ respectively, when the prior temperature is $T=0.01$:
\begin{equation}\nonumber
\left[
\begin{matrix}
0.5000 & 0.2000 & 0.0000 & 0.0000 & 0.0000 & 0.0000 & 0.0000 & 0.0806 & 0.2194 \\
0.0000 & 0.0000 & 0.0623 & 0.4476 & 0.0548 & 0.0000 & 0.1354 & 0.0806 & 0.2194 \\
0.0000 & 0.0000 & 0.0000 & 0.3952 & 0.0548 & 0.0000 & 0.1085 & 0.1952 & 0.2463 \\
0.0000 & 0.0000 & 0.0000 & 0.3429 & 0.0548 & 0.0000 & 0.0816 & 0.2476 & 0.2731 \\
0.0000 & 0.0000 & 0.0000 & 0.2905 & 0.0548 & 0.0000 & 0.0548 & 0.3000 & 0.3000 \\
\end{matrix}
 \right].
\end{equation}
Such evolution privileging shortest paths is also  apparent in Figure \ref{fig:boltzmann_T_0.01}. 
The numbers over the edges indicate the amount of mass that has traveled during $(t_i,t_{i+1})$ on that edge. For instance, we see that between $t=0$ and $t=1$, of the $0.5$ probability available in node $1$ at $t=0$, $0.448$ has moved to node $4$ while $0.052$ has moved to node $3$.
 In particular,  the mass accumulates on node 8 very quickly due to on the one hand to its proximity to nodes $1$ and $2$, where most of the initial mass is concentrated. On the other hand, node $8$ is close to the other ``sink" node $9$. Moreover, the mass hardly travels on longer paths: For instance, path $1-2-5-6-9$ {\em is not used at all}.

\begin{figure}[H]
    \begin{center}
    \includegraphics[width=15cm]{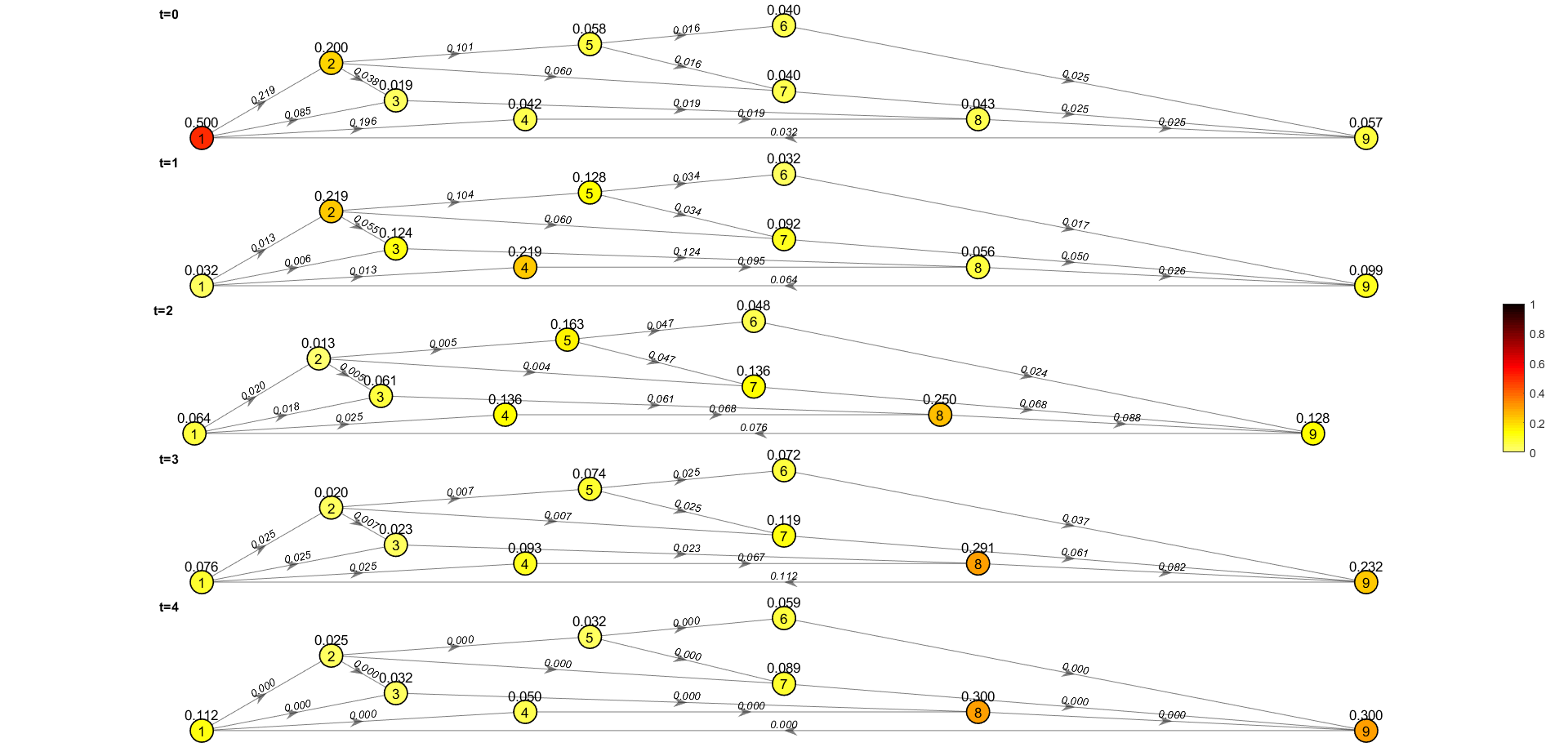}
    \end{center}
    \caption{Solution to the Schr\"odinger bridge problem with a Boltzmann prior T=100.}
    \label{fig:boltzmann_T_100}
\end{figure}

As a by-product of the optimal evolution, we get the following completion of the marginals at $t=0,4$:

\begin{eqnarray}\nonumber &&q_0^*(3)=q_0^*(4)=q_0^*(5)=q_0^*(6)=q_0^*(7)=0,\\\nonumber &&q_0^*(8)= 0.081, q_0^*(9)= 0.219\\\nonumber &&q_4^*(1)=q_4^*(2)=q_4^*(3)= 0, q_4^*(4)=0.29,\\\nonumber &&q_4^*(5)= 0.055, q_4^*(6)= 0, q_4^*(7)= 0.055.
\end{eqnarray}
 We now compare the above solution to the solution to same problem when the temperature in the Boltzmann prior has been raised to $T=100$, see the matrix below and Figure\ref{fig:boltzmann_T_100}:
 \begin{equation}\nonumber
\left[
\begin{matrix}
0.5000 & 0.2000 & 0.0187 & 0.0418 & 0.0581 & 0.0405 & 0.0405 & 0.0435 & 0.0571 \\
0.0323 & 0.2191 & 0.1237 & 0.2186 & 0.1279 & 0.0316 & 0.0920 & 0.0562 & 0.0986 \\
0.0637 & 0.0131 & 0.0613 & 0.1359 & 0.1634 & 0.0485 & 0.1361 & 0.2495 & 0.1284 \\
0.0756 & 0.0203 & 0.0229 & 0.0927 & 0.0742 & 0.0715 & 0.1191 & 0.2915 & 0.2323 \\
0.1119 & 0.0252 & 0.0320 & 0.0504 & 0.0317 & 0.0595 & 0.0894 & 0.3000 & 0.3000 \\
\end{matrix}
 \right],
\end{equation}

The high temperature makes so that the mass spreads as much as the topology of the graph allows, using also longer paths. For instance, some of the mass now reaches node 9 from node 1 along the path path $1-2-5-6-9$ which was not employed with $T=0.01$. Again, as a by-product, we read out a completion of the initial and final marginals $q_0^*$ and $q_4^*$, respectively.

\subsection{Ruelle-Bowen prior}

Next, we consider the example of the previous subsection with the addition of all self loops and an edge from node $9$ to node $1$ to make the graph strongly connected. We now have the graph topology:
\begin{figure}[H]
    \begin{center}
    \includegraphics[width=10cm]{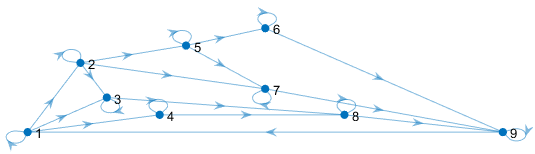}
    \end{center}
    \caption{Graph topology}
    \label{fig:topology_incomplete}
\end{figure}
\noindent
when the prior is given by the Ruelle-Bowen (RB) random walk measure \cite{Rue04,DL}, \cite[SectionIV]{chen2016robust}. This measure induces a uniform probability measure on paths of a fixed length between any two given nodes. Thus, the Ruelle-Bowen measure represents a natural choice for the prior in cases where we only know the topology of the graph. This distribution maximises the entropy rate for a random walker. The latter is bounded by the topological entropy rate which depends only on the topology of the graph and not its underlying probability distribution.
$$H_\mathcal{G}=\limsup_{t \to \infty}[\log|\{\text{paths of length t} \}/t$$

Let $A=(a_{ij})$ be the adjacency matrix of our graph $\mathcal{G}$. In this case the number of paths of length t is given by the nonzero entries of $A^t$. Let $\lambda_A$ be the spectral radius of A i.e the maximum modulus of the eigenvalues of A. We have that $H_\mathcal{G}=log(\lambda_A)$.

Let $u$,$v$ be the eigenvectors corresponding to $\lambda_A$ scaled so that their inner product is 1. We then define $$\nu_{RB}(i)=u_iv_i$$
which is invariant under the transition matrix
\[
R=(r_{ij}),\quad r_{ij}=\frac{v_j}{\lambda_A v_i}a_{ij}.
\]
We can now define the RB path measure as

\begin{equation}\label{RUBO}\mathcal {M}_{RB}(x_0,...,x_N)=\nu_{RB}(x_0)r_{x_0x_1}...r_{x_{N-1}x_N}.
\end{equation}

\begin{figure}
    \begin{center}
    \includegraphics[width=15cm]{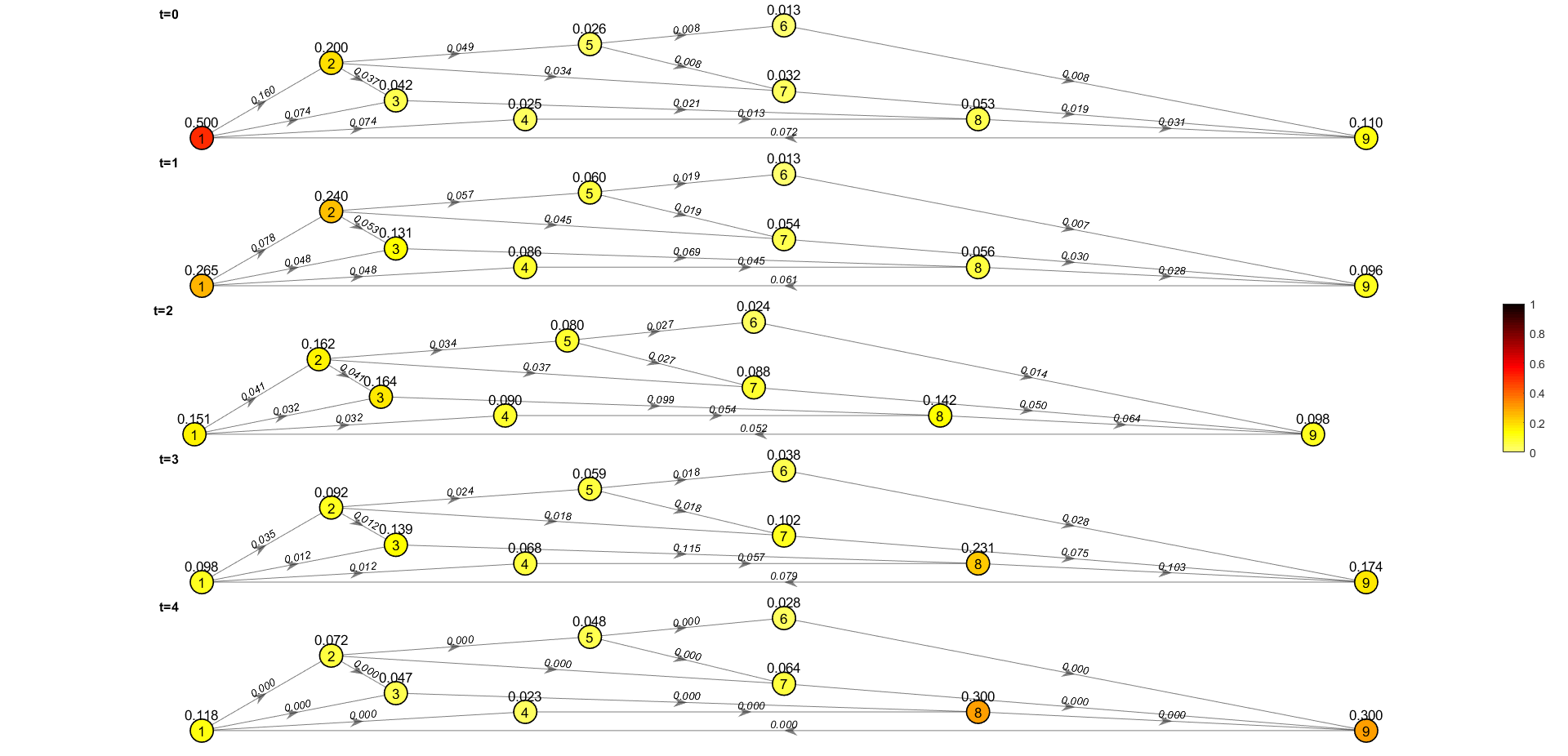}
    \end{center}
    \caption{Solution to the Schr\"odinger bridge problem with Ruelle-Bowen prior.}
    \label{fig:RB}
\end{figure}

Let us return to the  example of the previous subsection when the prior distribution is the {\em Ruelle-Bowen random walk measure}. Figure \ref{fig:RB} illustrates the distribution evolution. Similar amounts of mass travel on equal length paths between any two given nodes. 

For instance, there are $15$  paths of length $4$ from node $1$ to node $9$: 
\begin{eqnarray}\nonumber 1-2-5-6-9,\quad  1-1-2-7-9, \quad 1-2-2-7-9,\\\nonumber  1-2-5-7-9, \quad 1-2-7-7-9, \quad
1-1-3-8-9,\\\nonumber 1-2-3-8-9,\quad 1-3-3-8-9, \quad 1-1-4-8-9,\\\nonumber 1-4-4-8-9,\quad 1-3-8-8-9,\quad 1-4-8-8-9,\\\nonumber 1-2-7-9-9,\quad 1-3-8-9-9,\quad 1-4-8-9-9.
\end{eqnarray}
 Each of them has roughly (conditional) probability $0.09447$.

Next, we investigate what else can be said about the solution when the prior random vector happens to be a time window of a Markov chain.

\section{Markovian prior}\label{sec:prior}

Consider the situation where the prior distribution on paths $P$ is induced by a time-homogeneous Markov chain so that
\begin{equation}\label{prior}P(x_0,x_1,\ldots,x_N)=p_0(x_0)p_{x_0x_1}\cdots p_{x_{N-1}x_N},
\end{equation}
where $p_{ij}=\bbP(X_{t+1}=j|X_t=i)$.
Notice that both the Boltzmann distribution (\ref{Bolt}) and the Ruelle-Bowen measure (\ref{RUBO}) are of the (\ref{prior}) type. Then, in all SBP problems considered in this paper, it is interesting to observe that the solution $Q^*$ is also associated to  a Markov evolution. The new (time-varying) transition probabilities are as follows:

In the case when only (full or partial) information on the initial marginal is available $q^*_{ij}(t)= p_{ij}$.

In the case when only (full or partial) information on the final marginal is available
    \begin{equation}\label{opttran}q^*_{ij}(t)= \frac{\varphi(t+1,j)}{\varphi(t,i)}p_{ij},
    \end{equation}
    where 
    \[\varphi(t,i)=\frac{q^*_t(i)}{p_t(i)}, \quad t=0,1,\ldots,N.
\]
satisfies
\begin{equation}\label{phi1}\varphi(t,i)=\sum_{j}p(t,i;t+1,j)\varphi(t+1,j).
\end{equation}
The boundary condition at $t=N$ is given by $\nu_N(x_N)/p_N(x_N)$ is the full marginal knowledge case and by (\ref{endvarphi}) in the partial marginal knowledge case.

Consider now the full bridge problem with full information  on the two marginals $\nu_0$ and $\nu_N$. In such case, \cite{pavon2010discrete,GP}, the new transition probabilities are as in (\ref{opttran}) with $\varphi$ satisfying together with $\hat{\varphi}$ the Schr\"odinger system
\begin{eqnarray}
&&\varphi(t,i)=\sum_{j} p_{ij}\varphi(t+1,j),\\&&\hat{\varphi}(t+1,j)=\sum_{i}p_{ij}\hat{\varphi}(t,i).
\end{eqnarray}
for $t=0,1,\ldots,N-1$ with boundary conditions
\begin{equation}\label{bndconditions2}
\varphi(0,i)\cdot\hat{\varphi}(0,i)=\nu_0(i),\quad \varphi(N,j)\cdot\hat{\varphi}(N,j)=\nu_N(j),\quad \forall i, j\in\mathcal X.
\end{equation}
When only partial information of the marginals is available, the analysis of Subsections \ref{FBPI} and {\ref{ITAL} show that (\ref{bndconditions2}) need to be replaced with (\ref{BConestep}) supplemented with
\begin{eqnarray}\varphi(N,j)=1\; {\rm for}\; j\in\mathcal X_N^c\\\hat{\varphi}(0,i)=\frac{(1-\alpha)p_0(i)}{\sum_{i\in\mathcal X_0^c}p_0(i)\varphi(0,i)}\; {\rm for}\; i\in\mathcal X_0^c.
\end{eqnarray}

\section{Knowledge of the Moments of the Marginals}\label{moments}

We  consider the case where we  only some moments of the initial and/or final distribution are known. This represents the common case where these moments have been estimated from the data. Again, let $\mathcal P_{0N}$ be the family of probability distributions over $\mathcal X\times\mathcal X$. 


Suppose  first we only know $m_0$ and $m_N$ the means of the initial and final marginal, respectively. 
To avoid trivial cases, we assume  $1<m_0,m_N<n$. 
Using once more decomposition (\ref{DECO}), the maximum entropy problem becomes minimizing

\begin{equation}\bbD(q_{0N}\|p_{0N})=\sum_{x_0x_N}q_{0N}(x_0,x_N)\log \frac{q_{0N}(x_0,x_N)}{p_{0N}(x_0,x_N)}
\end{equation}

with respect to $q_{0N}\in\mathcal P_{0N}$ subject to the (linear) constraints

\begin{eqnarray}
\sum_{x_0,x_N}q_{0N}(x_0,x_N)x_0&=&m_0\\
\sum_{x_0,x_N}q_{0N}(x_0,x_N)x_N&=&m_N\\
\sum_{x_0x_N}q_{0N}(x_0,x_N)&=&1.\label{Normalization}
\end{eqnarray}

This now a standard,  strictly convex optimization problem \cite[Chapters 4, 5]{BV}  which can be solved through duality theory (strong duality holds). Indeed, the Lagrangian function has the form

\begin{eqnarray}\nonumber&&{\cal L}(q_{0N};\lambda,\mu,\theta)=\sum_{x_0x_N}q_{0N}(x_0,x_N)\log \frac{q_{0N}(x_0,x_N)}{p_{0N}(x_0,x_N)}\\\nonumber&&+\lambda\left[\sum_{x_0,x_N}q_{0N}(x_0,x_N)x_0-m_0\right]+\mu\left[\sum_{x_0,x_N}q_{0N}(x_0,x_N)x_N-m_N\right]\\&&+\theta\left[\sum_{x_0x_N}q_{0N}(x_0,x_N)-1\right]\label{lagrangian}
\end{eqnarray}
Setting the first variation equal to zero, we get the (sufficient) optimality condition:

\[q_{0N}^*(x_0,x_N)=\exp\left[-1-\theta\right]]\exp\left[-\lambda x_0\right]\exp\left[-\mu x_N\right]p_{0N}(x_0,x_N).\]
Multipliers $\theta$, $\lambda$ and $\mu$ can be determined solving the concave, unconstrained $3$-dimensional dual problem 
\[
{\rm maximize}_{(\theta,\lambda,\mu)\in(\bbR\times\bbR\times\bbR)}{\cal L}(q^*_{0N};\theta,\lambda,\mu) 
\]
through, e.g., gradient ascent. This case is actually so simple that one can argue directly about the existence and uniqueness of the multipliers studying a system of polynomial equations, see Appendix \ref{mean}. The case of knowledge of the first two moments of the initial and final marginals is outlined in Appendix \ref{variance}. Finally, half-bridges problems with only moment knowledge of one of the marginals can be solved in a similar fashion.

Consider the graph depicted in FIG. \ref{fig:topology_incomplete}.
We now consider the bridge problem with a Boltzmann prior (\ref{Bolt}) with $T=1$ in the case where only the means $m_0=1.5$ and $m_N=7$ of the initial and final marginal, respectively are known.
We observe that the optimal distribution is concentrated on minimum length paths due to the low $T$.

\begin{figure}[H]
    \begin{center}
    \includegraphics[width=15cm]{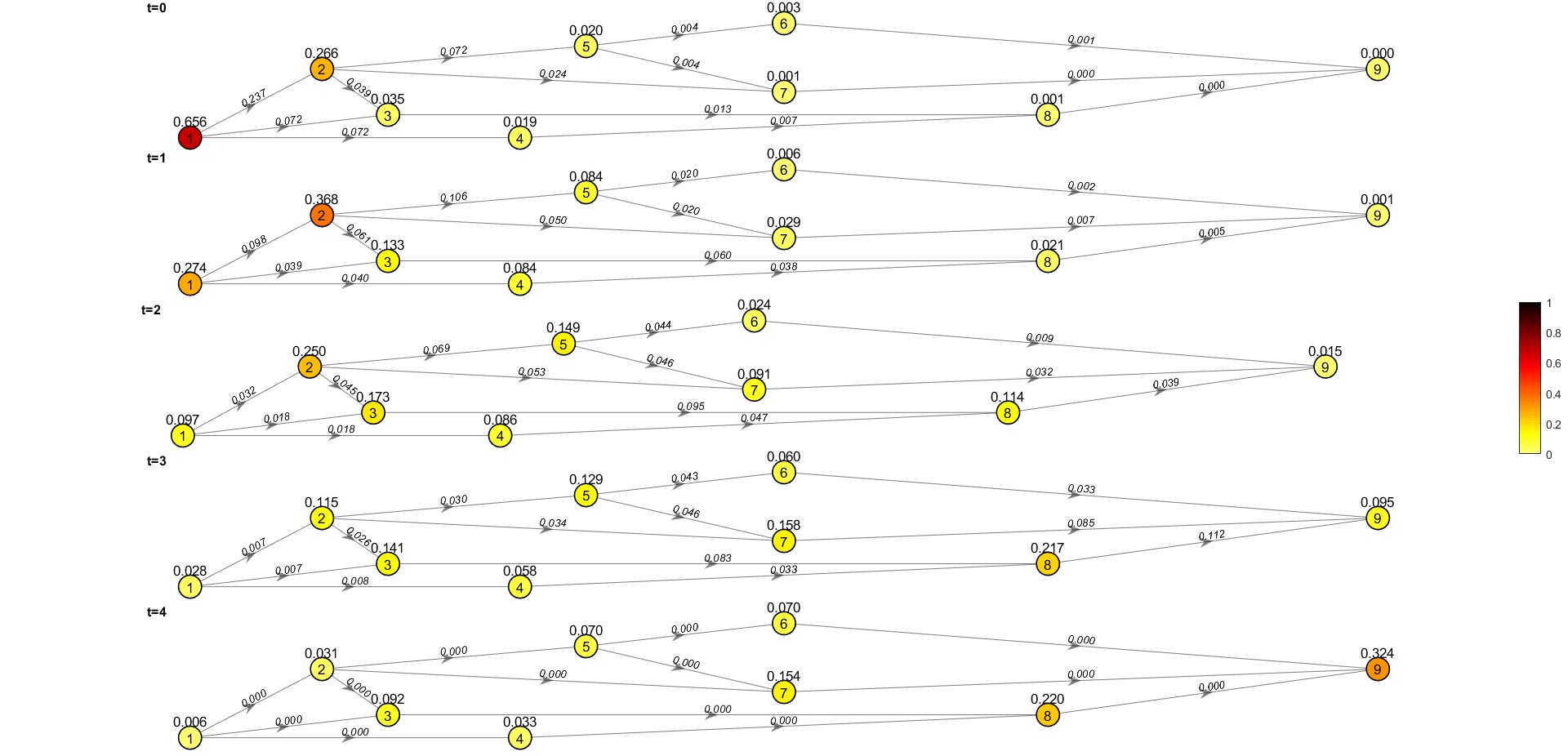}
    \end{center}
    \caption{Solution to a Schrodinger bridge problem with knowledge of the initial and final means}
    \label{fig:mean}
\end{figure}

\section{Conclusion}

We have shown that regularized optimal transport over networks can be solved also in the case of partial information about the initial and final distribution. Two cases were considered: In case one, one knows  one or both of the marginals only on a subset of the state space, see Section \ref{partial}. In the second case, only some moments of the initial and/or final marginal are known, see Section \ref{moments}.The solution in the first case can be computed through an efficient iterative scheme which represents a suitable modification of the classical Fortet-IPF-Sinkhorn algorithm. In the second case, standard convex optimization iterative methods can be used to solve the dual problem. Similar results in both cases can be established in the continuous-time continuous state space case, see Appendix \ref{continuous}. The  problem with moment information has been studied in the continuous setting without entropy regularization in \cite{Alfonsi2021}. In \cite{Elvander2020}, a multi-marginal regularized optimal transport problem is studied where the partial information on the marginals is available through various linear maps.
Albeit the mathematics is  admittedly simple, the potential of this result appears considerable given that the solution also entails a most likely completion of the initial and final marginals. 

\begin{acknowledgments}
We are very grateful to Esteban Tabak for proposing the problem considered in this paper. 
\end{acknowledgments}

\appendix 
\section{Mean of the initial and final marginal}\label{mean}

Using the constraints and the decomposition $p_{0N}(x_0,x_N) = p_0(x_0)p(0,x_0;N,x_N)$ we get

\begin{eqnarray}
\sum_{x_0}c\exp\left[-\lambda x_0\right]p_0(x_0)x_0\sum_{x_N}p(0,x_0;N,x_N)\exp\left[-\mu x_N\right]&=&m_0,\\
\sum_{x_N}c\exp\left[-\mu x_N\right]x_N\sum_{x_0}p(0,x_0;N,x_N)p_0(x_0)\exp\left[-\lambda x_0\right])&=&m_N,\\
\sum_{x_0x_N}c\exp\left[-\lambda x_0\right]\exp\left[-\mu x_N\right]p_{0N}(x_0,x_N)&=&1.
\end{eqnarray}

where $c = \exp\left[-1-\theta\right]$. We may now define 
\begin{eqnarray}
h(\mu,x_0) &:=& p_0(x_0)\sum_{x_N}p(0,x_0;N,x_N)\exp\left[-\mu x_N\right],\\
g(\lambda,x_N) &:=& \sum_{x_0}p(0,x_0;N,x_N)p_0(x_0)\exp\left[-\lambda x_0\right]).
\end{eqnarray}

Combining these definitions with the constraints, we get the following equations for $r:=exp(\lambda)$ and $s:=exp(\mu)$: 

\begin{eqnarray}
\sum_{x_0}h(\mu,x_0)x_0r^{-x_0}=\sum_{x_0}m_0h(\mu,x_0),\\
\sum_{x_N}g(\lambda,x_N)x_Ns^{-x_N}=\sum_{x_N}m_Ng(\lambda,x_N).
\end{eqnarray}

These, in turn can be transformed into the following polynomial equations for r and s.

\begin{eqnarray}
P(\mu)=\sum_{x_0}h(\mu,x_0)(x_0-m_0)r^{N-x_0}=0,\\
\hat{P}(\lambda)=\sum_{x_N}g(\lambda,x_N)(x_N-m_N)s^{N-x_N}=0.
\end{eqnarray}
Here the coefficients of the polynomials $P(\mu),\hat{P}(\lambda)$ are, respectively,
\[[h(\mu,1)(1-m_0), ... ,h(\mu,N)(N-m_0)], \quad [g(\lambda,1)(1-m_N), ... ,g(\lambda,N)(N-m_N]\]
Note that $h(\mu,x_0)$ and $g(\lambda,x_0)$ are positive for $x_0,x_N\in \mathcal X$. Hence, the sign of the coefficients is determined by the sign of $(i-m_0)$ and $(i-m_N)$, respectively. By Descartes Rule of signs $P(\mu)$ and $\hat{P}(\lambda)$ have exactly one positive real root.   We can find such roots through  an iterative algorithm. We define a root function $R_p$ for this positive real root. Let $R_P(\mu)$ and $R_{\hat{P}}(\lambda)$ be the root functions for the two polynomials. Then

\begin{eqnarray}
\mu  \overset{\cA}{\longrightarrow} \lambda, \quad \lambda &=& log(R_P(\mu))\\
\lambda \overset{\cB}{\longrightarrow} \mu,\quad \mu &=&log( R_{\hat{P}}(\lambda))
\end{eqnarray}
Establishing convergence of the above iteration would prove existence and uniqueness for the dual problem.

\section{Mean and variance of the initial and final marginal}\label{variance}

Let $m_{n,0}$ and $m_{n,N}$ be the nth moment of the initial and final marginal respectively. With information on the mean and variance we must minimize

\begin{equation}\bbD(q_{0N}\|p_{0N})=\sum_{x_0x_N}q_{0N}(x_0,x_N)\log \frac{q_{0N}(x_0,x_N)}{p_{0N}(x_0,x_N)}
\end{equation}

with respect to $q_{0N}\in\mathcal P_{0N}$ subject to the (linear) constraints

\begin{eqnarray}
\sum_{x_0,x_N}q_{0N}(x_0,x_N)x_0&=&m_{1,0}\\
\sum_{x_0,x_N}q_{0N}(x_0,x_N)x^2_0&=&m_{2,0}\\
\sum_{x_0,x_N}q_{0N}(x_0,x_N)x_N&=&m_{1,N}\\
\sum_{x_0,x_N}q_{0N}(x_0,x_N)x^2_N&=&m_{2,N}\\
\sum_{x_0x_N}q_{0N}(x_0,x_N)&=&1.\label{normalization}
\end{eqnarray}

A variational analysis similar to that of Section \ref{moments}, gives the following optimality condition:
\[q_{0N}^*(x_0,x_N)=p_{0N}(x_0,x_N)\exp\left[-1-\theta\right]]\exp\left[-\lambda x_0\right]\exp\left[-\mu x_N\right]\exp\left[-\alpha x^2_0\right]\exp\left[-\beta x^2_N\right]\]
We can then solve the unconstrained $5$-dimensional dual problem 
\[
{\rm maximize}_{(\theta,\lambda,\mu,\alpha,\beta)\in(\bbR)^4}{\cal L}(q^*_{0N};\theta,\lambda,\mu,\alpha,\beta) 
\]
through standard iterative methods.

\section{Continuous case with partial marginal information}\label{continuous}
The case of continuous time and state space can be dealt with in a similar fashion. Suppose we  consider probability measures induced on continuous functions $C([0,1];\bbR^n)$ by $n$-dimensional diffusion processes. Let $p_{01}(x,y)$ denote the joint initial-final density of the prior measure $P$. Suppose we only know the initial marginal density $\rho_0(x)$ on $\mathcal X_0$  and the final marginal density  $\rho_1(y)$ on $\mathcal X_1$. Then, a decomposition of relative entropy for path measures similar to (\ref{DECO}), see e.g. \cite[(4.6)]{chen2021stochastic}, reduces the problem to minimizing over joint densities $q_{01}$ on $\bbR^n\times\bbR^n$
\begin{equation}\label{IMSBP5}\bbD(q_{01}\|p_{01})=\int\int q_{01}(x,y)\log \frac{q_{01}(x,y)}{p_{01}(x,y)}dxdy
\end{equation}
subject to the linear constraints
\begin{eqnarray}\label{IMSBP6}\int q_{01}(x,y)dy&=&\rho_0(x),\quad x\in{\mathcal X_0},
\\
\label{IMSBP7}\int q_{01}(x,y)dx&=&\rho_1(y),\quad y\in{\mathcal X_1},\\\label{IMSBP8}
\int \int q_{01}(x,y)dxdy&=&1.
\end{eqnarray}
Introducing multipliers $\lambda(x),\mu(y),\theta$ for the three constraints, a variational analysis parallel to that of Subsection \ref{FBPI} and similar to \cite[Subsection 4.3]{chen2021stochastic} gives
\[
\frac{q_{01}^*(x,y)}{p(0,x;1,y)}=\left\{p_0(x)\exp\left[-1-\theta-\mathds{1}_{\mathcal X_0}(x)\lambda(x)\right]\right\}\left\{\exp\left[-\mathds{1}_{\mathcal X_1}~(y)\mu(y)\right]\right\}=\hat{\varphi}(x)\cdot\varphi(y).
\]
Here
\begin{equation}\nonumber\hat{\varphi}(x)= \left\{\begin{array}{ll} p_0(x)\exp [-1-\theta-\lambda(x)],& {\rm if}\;x\in \mathcal X_0\\ p_0(x_0)\exp [-1-\theta],& {\rm if}\;x\in (\mathcal X_0)^c\end{array}\right.,
\end{equation}
and
\begin{equation}\nonumber\varphi(y)= \left\{\begin{array}{ll} \exp \left[-\mu(y)\right],& {\rm if}\;y\in \mathcal X_1\\ 1,& {\rm if}\;y\in (\mathcal X_1)^c\end{array}\right.
\end{equation}
The iterative algorithm of Subsection \ref{ITAL} is replace by
Consider the following four maps
\begin{eqnarray}\hat\varphi(0,x)  \overset{\cE^\dagger}{\longrightarrow} \hat\varphi(1,y),\quad \hat\varphi(1,y) &=& \int p(0,x;1,x_y)\hat{\varphi}(0,x_)dx\\
 \hat\varphi(1,y)\overset{\cD_{1}}{\longrightarrow}\varphi(1,y),\quad \varphi(1,y)&=& \left\{\begin{array}{ll} \frac{\rho_1(y)}{\hat{\varphi}(1,y)},& {\rm if}\;y\in \mathcal X_1\\ 1,& {\rm if}\;y\in (\mathcal X_1)^c\end{array}\right.\\
\varphi(1,y)
\overset{\cE}{\longrightarrow} \varphi(0,x),\quad \varphi (0,x)&=&\int p(0,x_0;1,y)\varphi(1,y)\\
 \varphi(0,x)\overset{\cD_{0}}{\longrightarrow}\left(\hat\varphi(0,x)\right)_{\rm next},\quad\hat{\varphi}(0,x)^{{\rm next}}&=& \left\{\begin{array}{ll} \frac{\rho_0(x)}{\varphi(0,x)},& {\rm if}\;x\in \mathcal X_0\\ c_0p_0(x_0),& {\rm if}\;x_0\in (\mathcal X_0)^c\end{array}\right.\label{C8}
\end{eqnarray}
Here, $c_0=\exp[-1-\theta]$ must be such that
\[
\int_{\bbR^n}\hat\varphi(0,x)\cdot \varphi(0,x)dx=1.
\]
Let 
\[
\int_{\mathcal X_0}\rho_0(x)dx=\alpha, \quad 0<\alpha<1. 
\]
It follows from (\ref{C8}) that $c_0$ is given by
\[c_0=\frac{1-\alpha}{\int_{\mathcal X_0}p_0(x)\varphi(0,x)dx}.
\]
Hence, map (\ref{C8}) now reads
\begin{equation} \varphi(0,x)\overset{\cD_{0}}{\longrightarrow}\left(\hat\varphi(0,x)\right)_{\rm next},\quad\hat{\varphi}(0,x)^{{\rm next}} = \left\{\begin{array}{ll} \frac{\rho_0(x)}{\varphi(0,x)},& {\rm if}\;x\in \mathcal X_0\\ \frac{(1-\alpha)p_0(x)}{\int_{\mathcal X_0}p_0(x)\varphi(0,x)dx},& {\rm if}\;x_0\in (\mathcal X_0)^c\end{array}\right.
\end{equation}
Convergence of the composition in the Hilbert metric requires some further work and care due to some difficulties intrinsic of the continuous case, see \cite{CheGeoPav15a}. In alternative, one can prove convergence of Fortet's original iterative system, see \cite{fortet1940,EP,leo2019}.

\begin{thebibliography}{99}
\bibitem{Alfonsi2021} {\sc A. Alfonsi, C. Rafa{\"e}l, V. Ehrlacher and D. Lombardi} Approximation of optimal transport problems with marginal moments constraints, {\em Mathematics of Computation}, {\bf 90} n.328, (2021), 689--737.

\bibitem{BCCNP} {\sc J. Benamou, G. Carlier, M. Cuturi, L. Nenna and G. Peyr\'e}, 
Iterative Bregman projections for regularized transportation problems,  {\em SIAM Journal on Scientific Computing}, {\bf 37}, n.2, (2015), A1111--A1138.

\bibitem{BV} S. Boyd and L. Vandenberghe, {\em Convex Optimization}, Cambridge University Press, 2004.
\bibitem{CheGeoPav15a} {\sc Y. Chen, T. T. Georgiou and M. Pavon}, Entropic and displacement interpolation: a computational approach using the {H}ilbert metric, {\em SIAM Journal on Applied Mathematics}, {\bf 76}, n.6, (2016), 2375--2396. 


	
\bibitem{chen2016robust} {\sc Y. Chen, T. T. Georgiou, M. Pavon and A. Tannenbaum}, Robust transport over networks, {\em IEEE Transactions on Automatic Control}, {\bf 62}, n. 9, (2016), 4675--4682.
	

\bibitem{chen2017efficient} {\sc Y. Chen, T. T. Georgiou, M. Pavon and A. Tannenbaum}, Efficient robust routing for single commodity network flows, {\em IEEE Transactions on Automatic Control}, {\bf 63}, n.7, (2017), 2287--2294.

\bibitem{chen2021stochastic}Y.\ Chen, T.T.\ Georgiou and M.\ Pavon, Stochastic control liaisons: Richard Sinkhorn meets Gaspard Monge on  a Schr\"{o}dinger bridge, {\em SIAM Review}, {\bf 63}, n.2, (2021), 249--313. 

\bibitem{chen2021most} {\sc Y.Chen, Yongxin, T.T. Georgiou,  and M. Pavon}, The most likely evolution of diffusing and vanishing particles: Schr\"{o}dinger Bridges with unbalanced marginals, {\em SIAM J. Contr. Optimiz.}, {\bf 60}, n. 4, (2022), 2016--2039.

	

\bibitem{chizat2018scaling} {\sc L.Chizat, G. Peyr{\'e}, B. Schmitzer and F.-V. Vialard}, Scaling algorithms for unbalanced optimal transport problems, {\em Mathematics of Computation}, {\bf 87}, n. 314, (2018), 2563-2609.

\bibitem{CGKS} {\sc S. Courtain, G. Guex, I. Kivimäki and M. Saerens},   Relative entropy-regularized optimal transport on a graph: a new algorithm and an experimental comparison, {\em International Journal of Machine Learning and Cybernetics}, {\bf 14}(4), (2023), 1365-1390.

\bibitem{CLZ}{\sc J. Cui, S. Liu, and  H. Zhou}, What is a stochastic Hamiltonian process on finite graph? An optimal transport answer, {\em Journal of Differential Equations}, {\bf 305}, (2021), 428--457.
	
\bibitem{DL} J.-C. Delvenne and A.-S. Libert, Centrality measures and thermodynamic
formalism for complex networks, {\em Physical Review E}, {\bf 83},
no. 4, (2011) p. 046117.	
  
\bibitem{demingstephan1940}
{\sc W.~E. Deming and F.~F. Stephan}, On a least squares adjustment of a sampled frequency table when the expected marginal totals are known, {\em Ann. Math. Statist.}, {\bf 11} (1940), 427--444.

\bibitem{DB} {\sc X. Duan and F. Bullo},  Markov chain–based stochastic strategies for robotic surveillance. {\em Annual Review of Control, Robotics, and Autonomous Systems}, n. 4,(2021), 243--264.

\bibitem{Elvander2020}, {\sc F. Elvander, I. Haasler, A. Jakobsson, J. Karlsson}, Multi-marginal optimal transport using partial information with applications in robust localization and sensor fusion, {\em Signal Processing}, {\bf 171}, (2020), 107474.

\bibitem{EP} {\sc M. Essid and M. Pavon}, Traversing the Schr\"odinger Bridge strait: Robert Fortet's marvelous proof redux, {\em J. Optim. Theory and Applic.}, {\bf 181}, n.1, (2019), 23--60.

\bibitem{fortet1940}
{\sc R.~Fortet}, {\em R\'esolution d'un syst\`eme d'equations de {M}.
  {S}chr\"{o}dinger}, J. Math. Pure Appl., IX (1940), ~83--105.

  
\bibitem{GP} {\sc T. T. Georgiou and M. Pavon}, Positive contraction mappings for classical and quantum Schr\"{o}dinger systems, {\em J. Math. Phys.}, {\bf 56}, 033301 (2015).

\bibitem{HCK} {\sc I. Haasler, Y. Chen and J. Karlsson},  Optimal steering of ensembles with origin-destination constraints. {\em IEEE Control Systems Letters}, {\bf 5}, n.3, (2020), 881--886.

\bibitem{HRCK} {\sc I. Haasler, A. Ringh, Y. Chen, Y. and J. Karlsson}, Multimarginal optimal transport with a tree-structured cost and the schrodinger bridge problem, {\em SIAM J. Control and Optimization}, {\bf 59},n.4, (2021), 2428--2453.

\bibitem{KDO} {\sc P. Koehl, M. Delarue, and H. Orland}, Physics approach to the variable-mass optimal-transport problem, {\em Physical Review E}, {\bf 103},n.1,(2021), 012113.

\bibitem{LCFS} {\sc P. Leleux, S. Courtain, K. Françoisse and M. Saerens}, Design of biased random walks on a graph with application to collaborative recommendation. {\em Physica A: Statistical Mechanics and its Applications}, {\bf 590}, (2022), 126752.



\bibitem{leoD} {\sc C. L\'eonard}, Lazy random walks and optimal transport on graphs, {\em The Annals of Probability}, {\bf 44}, n.3, (2016) 1864--1915.

\bibitem{leo2019} C. L\'eonard, Revisiting Fortet's proof of existence of a solution to the Schr\"odinger system, ArXiv e-prints, arXiv:1904.13211.

\bibitem{liero2018optimal} {\sc M. Liero, A. Mielke and G. Savar{\'e}}, Optimal entropy-transport problems and a new {H}ellinger--{K}antorovich distance between positive measures, {\em Inventiones mathematicae}, {\bf  211}, n. 3, (2018), 969--1117.
	
\bibitem{Mascherpa2023} {\sc M. Mascherpa, I. Haasler, B. Ahlgren, J. Karlsson}, Estimationg pollution spread in water networks as a Schr\"odinger bridge problem with partial information, {\em European Journal of Control}, {\bf 74} (2023),100846.

\bibitem{pavon2010discrete}
{\sc M.~Pavon and F.~Ticozzi}, Discrete-time classical and quantum
  {M}arkovian evolutions: Maximum entropy problems on path space, {\em J. 
  Mathematical Physics}, {\bf 51} (2010), ~042104.

  \bibitem{Pey19}{\sc G. Peyr\'e and M. Cuturi}, Computational Optimal Transport, {\em Foundations and Trends in Machine Learning}, {\bf 11}, n. 5-6, (2019), 1--257.
	
\bibitem{Rue04} {\sc D.Ruelle}, {\em Thermodynamic formalism: The mathematical structure of equilibrium statistical mechanics}, (2004), Cambridge University Press.


\bibitem{sanov1957largedeviations}
{\sc I.~S. Sanov}, {\em On the probability of large deviations of random
  magnitudes (in russian)}, Mat. Sb. N. S., 42 (84) (1957), ~11--44.

\bibitem{S1} {\sc E. Schr\"{o}dinger}, \"{U}ber die Umkehrung der Naturgesetze,
{\em Sitzungsberichte
der Preuss Akad. Wissen. Berlin, Phys. Math. Klasse}, {\bf 10}, (1931), 144--153.

\bibitem{S2} {\sc E. Schr\"{o}dinger}, Sur la th\'{e}orie relativiste de
l'\'{e}lectron et l'interpretation de
la m\'{e}canique quantique, {\em Ann. Inst. H. Poincar\'{e}} {\bf 2}, (1932),
269.


\bibitem{Sinkhorn1964}
{\sc R.~Sinkhorn}, {\em A relationship between arbitrary positive matrices and
  doubly stochastic matrices}, The Annals of Mathematical Statistics, 35
  (1964), 876--879.

\bibitem{WFBR} {\sc A.Waqas, H. Farooq, N. C. Bouaynaya and G. Rasool},  Exploring robust architectures for deep artificial neural networks, {\em Communications Engineering}, {\bf 1}, n.1, (2022), 46.

\bibitem{Z} {\sc H. Zhou}, Optimal transport on networks, {\em IEEE Control Systems Magazine}, {\bf 41} n. 4, (2021), 70--81.


\end{thebibliography}
\end{document}